\documentclass[a4paper, 12pt]{amsart}
\usepackage{amsmath, amssymb, eucal, amscd, amstext, enumerate}
\usepackage{amsfonts,amsthm}
\usepackage{mathrsfs}

\theoremstyle{plain}
\newtheorem{thm}{Theorem}[section]
\newtheorem{lem}[thm]{Lemma}

\newtheorem{prop}[thm]{Proposition}

\newtheorem{defn}[thm]{Definition}
\newtheorem{rem}[thm]{Remark}
\newtheorem{rem-ntn}[thm]{Remark and Notation}

\newenvironment{prf}{{\noindent \textbf{Proof:}\ }}{\hfill $\Box$\\ \smallskip}

\numberwithin{equation}{section}

\newcommand{\smnoind}{{\smallskip\noindent}}

\newcommand{\id}{{\rm id}}

\newcommand{\ad}{{\rm ad}}
\newcommand{\Pol}{{\rm Pol}}

\newcommand{\abs}[1]{\left\vert#1\right\vert}

\newcommand{\CL}{\mathcal{L}}

\newcommand{\CR}{\mathcal{R}}

\newcommand{\CB}{\mathcal{B}}

\newcommand{\KH}{\mathfrak{H}}

\newcommand{\BN}{\mathbb{N}}

\newcommand{\BC}{\mathbb{C}}
\newcommand{\BF}{\mathbb{F}}

\begin{document}

\title[A Smiley-type theorem for spectral operators of finite type]{A Smiley-type theorem for spectral operators of finite type}

\author{Xiao Chen}
\address{School of Mathematics and Statistics, Shandong University, Weihai 264209, P.R. China}
\email{chenxiao@sdu.edu.cn}
\thanks{Xiao Chen is supported by the NSF of China (Grants 11701327).}

\author{Jian-Jian Jiang}
\address{School of Mathematics and Physics, Ningde Normal University, Ningde 352100, P.R. China}
\email{j.j.jiang@foxmail.com}
\thanks{Jian-Jian Jiang is supported by the NSF of China (Grants 11801290).}

\author{Xiaolin Li$^\dag$}
\address{School of Mathematics and Statistics, Shandong University, Weihai 264209, P.R. China}
\email{lixiaolin26@163.com}
\thanks{}




\begin{abstract}
In this short article, we mainly prove that, for any spectral operator $A$ of type $m$ on a complex Hilbert space, if a bounded operator $B$ lies in the collection of  bounded linear operators that are in the $k$-centralizer of every bounded linear operator in the $l$-centralizer of $A$, where $k\leqslant l$ is two arbitrary positive integers satisfying $l\geqslant k$ as well as $l\geqslant 2m+1$, then $B$ must belong to the von Neumann algebra generated by $A$ and identity operator. 
This result generalizes a matrix commutator theorem proved by M.\ F.\ Smiley. 
For this aim, Smiley-type operators are defined and studied.  
\vspace{04pt}

\vspace{04pt}

\end{abstract}

\begin{thanks} {$\dag$ the corresponding author: lixiaolin26@163.com}
\end{thanks}

\maketitle

\section{Introduction and preliminaries}\label{sec:intro-pre}
\smallskip

In 1960, M.\ F.\ Smiley in \cite{Smi1961} proved the following important and interesting theorem.
\emph{\begin{quote}
{\bf Classical Smiley's Theorem:}
Let $\mathbb{F}$ be an algebraically closed field and the characteristic is 0 or at least $n$, and let $M_n(\mathbb{F})$ be the ring of all $n$ by $n$ matrices with elements in $\mathbb{F}$. Let $A,B\in M_n(\mathbb{F})$ be such that for some positive integer $s$, $\ad_{A}^{s}(X)=0$ for $X$ in $M_n(\mathbb{F})$ implies that $\ad_{X}^{s}(B)=0$. Then $B$ is a polynomial in $A$ with coefficients in $\mathbb{F}$. Here $\ad$ will be defined below.
\end{quote}}

After Smiley, in \cite{Rob1964}, D.\ W.\ Robinson proved the above theorem is also valid in case $\BF$ is not algebraically closed, and had given a final and complete form of Smiley's theorem for matrix algebras(cf. \cite{Rob1964} or \cite[pp.114-115]{Rob1976}). In this paper, we will be devoted to generalize Smiley's theorem to infinite dimensional complex Hilbert space.

Let $\KH$ be a Hilbert space over complex number field $\BC$, and $B(\KH)$ be the $C^*$-algebra of all bounded operators on $\KH$. On $B(\KH)$, we can define a Lie product $[X,Y]:=XY-YX$ for any two operators $X,\ Y\in B(\KH)$. So $B(\KH)$ can be viewed as a Lie algebra. For any $T\in B(\KH)$, the adjoint operator of $T$ is denoted by $T^*$. The algebra generated by a subset $\mathcal{S}\subseteq B(\KH)$, denoted by $\langle\mathcal{S}\rangle$, means the smallest algebra containing $\mathcal{S}$. 

Denote by $\sigma(A)$ the spectrum of $A\in B(\KH)$.

For any operator $Z\in B(\KH)$, we can define the corresponding {\bf left (respectively, right) multiplier} on $B(\KH)$ by 
$$\CL_{Z}(X):=ZX \ (\text{\rm respectively, } \CR_{Z}(X):=XZ),\text{\rm for all }X\in B(\KH).$$

For any $A\in B(\KH)$, the map $\ad_A: B(\KH)\rightarrow B(\KH)$ determined by the formula 
$$\ad_A(X):= [A,X],\ \forall X\in B(\KH)$$ is a linear operator on $B(\KH)$. 
The map which to every operator $A$ assigns the operator
$\ad_A$ is called the {\bf adjoint representation} of $B(\KH)$. 
For any positive integer $s$, we denote by $\ad^s_A$ the $s$-multiple composition of the operator $\ad_A$,   
that is, $$\ad^s_A(X):=[A,[A,[A,[A,\cdots[A,X]\cdots]]]],$$ where $[A,\cdot]$ is repeated $s$ times.

Denote by $\BN$ the set of natural numbers. and by $\BN^+$ the set $\BN\setminus\{0\}$.

If $\CB$ is a subset of $B(\KH)$, we denote by $\CB'$ the {\bf commutant} of $\CB$, i.e. $\CB':=\{Y\in B(\KH)\ |\ XY=YX, \forall X\in\CB\}$. The {\bf double commutant} $\CB''$ of $\CB$ is $(\CB')'$. 
Similarly, the {\bf $3$-commutant} $\CB'''$ of $\CB$ is $(\CB'')'$, and the {\bf $s$-commutant} $\CB^s$ of $\CB$ is $(\CB^{s-1})'$ for any positive integer $s\geqslant 4$.  
Also, for any $s\in\BN^+$, we define {\bf $s$-centralizer} of $\CB$ by $C_s(\CB)=\{X\in B(\KH)\ |\ \ad_Y^s(X)=0, \forall Y\in\CB\}$. 
These concepts also appeared in some more early references, such as \cite[Chapter 4]{Mur1990} and \cite[pp. 113]{Rob1976}.
In particular, when $\CB$ is a singleton set $\{Z\}$ in $B(\KH)$, we can abbreviate $C_s(\{Z\})$ as $C_s(Z)$. Note that $C_1(Z)=C'(Z)$.

For any operator $A$ in $B(\KH)$, we denote by $\Pol(A)$ the ring of all polynomials in $A$ with coefficients in $\BC$.

\begin{defn}\label{defn:smiley-type-oper}
A operator $A\in B(\KH)$ is called {\bf $(k,l)$-type Smiley operator} if there exist two positive integers $k,l\in\BN^+$ such that $C_k(C_l(A))\subseteq VN(A)$ where $VN(A)$ is the von Neumann algebra in $B(\KH)$ generated by $A$ and the identity operator $\id_\KH$ on $\KH$. In addition, a $(k,l)$-type Smiley operator is said to be {\bf proper} if $C_k(C_l(A))$ is contained in the subalgebra $\Pol(A)$ of $VN(A)$.
\end{defn}
If $\KH_n$ is a $n$-dimensional complex Hilbert space, which is isomorphic to $\BC^n$, for an arbitrary positive integer $n$, then it is clear that $B(\KH_n)$ equals to the matrix Lie algebra $M_n(\BC)$.  For $B(\KH_n)$, using the terminologies above, we can rewrite the classical Smiley's theorem over a special field $\BC$ as follows (cf. \cite[pp.113--115]{Rob1976}):
\emph{\begin{quote}
For any positive integer $s$ and $A\in M_n(\BC)$, one has $C_s(C_s(A))\subseteq\Pol(A)$.
\end{quote}}
In other words,  \emph{every $n$ by $n$ matrix over $\BC$ is a proper $(s,s)$-type Smiley operator on $\BC^n$ for all positive integer $s\in\BN^+$}.  

More interestingly, von Neumann's \emph{double commutant theorem} (see \cite[Theorem 4.1.5]{Mur1990}) actually tells us that every $A\in B(\KH)$ is $(1,1)$-type Smiley operator. Then the following question arises naturally.
\emph{\begin{quote}
Which operators on a Hilbert space are $(k,l)$-type Smiley operators for some given positive integers $k,l\in\BN^+$?
\end{quote}}
In the present article, we will  partially answer this question for spectral operators on a complex Hilbert space. 

Loosely speaking, a spectral operator is an operator which has a spectral reduction, that is, it can be reduced by a family of spectral projections. These projection is also known as the resolution of the identity or the spectral resolution of the given operator. It has been observed in \cite{Dun1952} that the spectral reduction is just the Jordan canonical form in classical matrix theory. In other words, every matrix is a spectral operator. Besides, another famous example of a spectral operator is a normal Hilbert space operator which has spectral measures and spectral resolution (cf. \cite[Section 4.3]{Sol2018}).  Here we use an equivalent formulation of spectral operators.

\begin{defn}\cite[Theorem 5, Section 4, Chapter XV]{Dun-Sch1988}\label{defn:spec-oper}
An operator $T\in B(\KH)$ is called {\bf spectral operator} if it has a unique decomposition $T=S+N$ of a bounded {\bf scalar type operator} $S$ and a {\bf quasi-nilpotent operator} $N$ commuting with $S$. Here a scalar type operator $S$ means $S$ has a unique spectral resolution $E$ for which $S=\int_{\sigma(S)}z\, {\rm d}E(z)$, and a quasi-nilpotent operator $N$ is an operator with a zero spectrum. Note that $T$ and $S$ has the same spectrum and the same spectral resolution.
\end{defn}
For more details about spectral operators, the readers can refer to the good survey \cite{Dun1958} and the famous book \cite{Dun-Sch1988}.

Clearly, the most common quasi-nilpotent operator is the {\bf nilpotent operator}, the $k$th-power of which is the zero operator for some positive integer $k$. 
Moreover, in \cite{Wer1954}, J.\ Wermer had shown that \emph{the scalar type operators in a Hilbert space are those operators which are similar to normal operators}. 
Hence, we mainly consider a special and big family of spectral operators in the following definition.

\begin{defn}\cite[Definition 3, Section 5, Chapert XV]{Dun-Sch1988}\label{defn:spec-oper-fin} An operator $T\in B(\KH)$ is called {\bf a spectral operator of type $m$}, if the nilpotent operator $N$ in the {\bf canonical decomposition} $T=S+N$ in Definition~\ref{defn:spec-oper} is of {\bf order $m+1$} for some $m\in\BN$, that is, the number $m+1$ is the smallest positive integer such that $N^{m+1}=0$. We denote by $\mathscr{S}_m(\KH)$ the set of all bounded spectral operators of type $m$ on $\KH$, and by $\mathscr{S}_{\rm fin}(\KH)$ the set of all bounded spectral operators of finite type on $\KH$, i.e. $\mathscr{S}_{\rm fin}(\KH)=\bigcup_{m\in\BN}\mathscr{S}_m(\KH)$.
\end{defn}

It can be apparently seen that the operators in $\mathscr{S}_{\rm fin}(\KH)$ is the most direct generalization of classical matrices to the infinite dimensional case, since the Jordan canonical form of a matrix is a sum of a diagonal matrix and a nilpotent matrix. This is another reason why we here mainly consider the subclass $\mathscr{S}_{\rm fin}(\KH)$ of spectral operators in $B(\KH)$.

The following theorem is the main result of our paper. It can be viewed as a generalization of Smiley's theorem onto $\mathscr{S}_{\rm fin}(\KH)$. 

\begin{thm}\label{thm:spec-oper-smiley}
For any $m\in\BN$, every operator in $\mathscr{S}_m(\KH)$ is a $(1,1)$-type Smiley operator, and also  a $(s,s)$-type Smiley operator for every positive integer $s\geqslant 2m+1$.
\end{thm}

Since both $C_k(C_l(A))\subseteq C_l(C_l(A))$  and $C_k(C_l(A))\subseteq C_k(C_k(A))$ hold for every element $A\in B(\KH)$ and any two positive integers $k,l$ satisfying $l\geqslant k$, a more general result, as a direct consequence of Theorem~\ref{thm:spec-oper-smiley-kl}, is immediately obtained.

\begin{thm}\label{thm:spec-oper-smiley-kl}
For any $m\in\BN$, every operator in $\mathscr{S}_m(\KH)$ is a $(1,1)$-type Smiley operator, and also a $(k,l)$-type Smiley operator for any two positive integers $k$ and $l$ satisfying $l\geqslant k$ as well as $l\geqslant 2m+1$.
\end{thm}

In Section~\ref{sec:prf}, we prove two essential lemmas, and apply them to prove Theorem~\ref{thm:spec-oper-smiley}. In Section~\ref{sec:ex}, we give an example of Smiley operator, which is provided by a kind of non-spectral operator.

\medskip
\section{Proof of the main theorem}\label{sec:prf}
\smallskip

Before proving our main theorem, we need to list some vital lemmas. 

The proof of the first lemma now follows from similar lines of argument as that of \cite[Lemma 1]{Smi1961}.
For completeness, we present the argument here.

\begin{lem}\label{lem:s-1}
If $A$ is a bounded normal operator on a Hilbert space $\KH$, then, for any $X\in B(\KH)$, we have $\ad_{A}^{s}(X)=0$ for some positive integer $s$ implies that  $\ad_{A}(X)=0$.
\end{lem}

\begin{prf}
Since $A\in B(\KH)$ is a normal operator, there exists a unique spectral resolution $E$ such that $A=\int_{\sigma(A)} z\, {\rm d}E(z)$. 
Due to $[A,E(\Omega)]=0$ for all Borel sets $\Omega$ of $\sigma(A)$, the Jacobi identity for Lie product shows that $\ad_{A}\ad_{E(\Omega)}(X)=\ad_{E(\Omega)}\ad_{A}(X)$ for all $X$ on $B(\KH)$. 
So, since it is well known that $T\in B(\KH)$ commutes with both $A$ and $A^*$ if and only if $TE(\Omega)=E(\Omega)T$ for all Borel sets $\Omega$ of $\sigma(A)$ (cf. \cite[Chapter II, Theorem 2.5.5]{Mur1990}),  it follows from the classical Fuglede's theorem \cite[Chapter IV, Theorem 4.10]{Sol2018} that $\ad_{A}^{s}(X)=[\ad_{A}^{s-1}(X),A]=0$ gives $$\ad_{A}^{s-1}\ad_{E(\Omega)}(X)=\ad_{E(\Omega)}\ad_{A}^{s-1}(X)=[\ad_{A}^{s-1}(X),E(\Omega)]=0.$$ 
By induction, we have $\ad_{A}\ad_{E(\Omega)}(X)=0$ from which $\ad_{E(\Omega)}^{2}(X)=0$ follows at once. But $\ad_{E(\Omega)}^{2}(X)=E(\Omega)(E(\Omega)X-XE(\Omega))-(E(\Omega)X-XE(\Omega))E(\Omega)=0$ yields $E(\Omega)X=XE(\Omega)$ upon right and left multiplication by $E(\Omega)$. Thus $\ad_{E(\Omega)}(X)=0$ for all spectral projections $E(\Omega)$ and consequently $\ad_{A}(X)=0$.
This completes the proof.
\end{prf}

\begin{lem}\label{lem:im-real-commut}
Let $A\in B(\KH)$ be a self-adjoint operator. For any $X\in B(\KH)$, one has that $\ad_A(X)=0$ implies both of $\ad_A(\Re X)=0$ and $\ad_A(\Im X)=0$, where the real part $\Re X=\frac{X+X^{\ast}}{2}$ and the imaginary part $\Im X=\frac{X-X^{\ast}}{2i}$ are both self-adjoint operators.
\end{lem}

\begin{prf}  Before proving this lemma, we first list two well-known facts.

\smnoind
{\bf Fact 1:} \emph{if $A$ and $B$ are both self-adjoint operators on $\KH$, then $(\ad_{A}(B))^{\ast}=-\ad_{A}(B)$.} 
Indeed,  it can be easily checked that $(\ad_{A}(B))^{\ast}=(AB-BA)^*=B^{\ast}A^{\ast}-A^{\ast}B^{\ast}=BA-AB=-\ad_{A}(B)$.

\smnoind
{\bf Fact 2:} \emph{if $A\in B(\KH)$ and $A^{\ast}=-A$, then there exists a unique self-adjoint operator $B$ such that $A=iB$, and $B$ is just $\Im A$.}
For clarifying {\bf Fact 2}, we set $A_{1}=\Re A$ and $\Im A=A_{2}$. Then $A^{\ast}=A_{1}-iA_{2}$. So $A+A^{\ast}=A_{1}+iA_{2}+A_{1}-iA_{2}=2A_{1}=0$, i.e., $A_{1}=0$. Hence $A=iA_{2}$. Uniqueness is obvious. 

Now, let us come to prove this lemma. Because $A$ is a self-adjoint operator, by {\bf Fact 1}, we know that 
$(\ad_{A}(\Re X))^{\ast}=-\ad_{A}(\Re X)$ and $(\ad_{A}(\Im X))^{\ast}=-\ad_{A}(\Im X)$. 
Hence, by {\bf Fact 2}, there exists two self-adjoint operators $Y_R$ and $Y_I$ such that $\ad_{A}(\Re X)=iY_R$ and $\ad_{A}(\Im X)=iY_I$. Then $\ad_A(X)=\ad_A(\Re X)+i\ad_A(\Im X)=iY_R-Y_I$, which means that $\Re (\ad_A(X))=-Y_I$ and $\Im (\ad_A(X))=Y_R$. This shows that, if $\ad_A(X)=0$, then $\ad_{A}(\Re X)=iY_R=i\Im (\ad_A(X))=0$ and $\ad_{A}(\Im X)=iY_I=-i\Re (\ad_A(X))=0$. This lemma is obtained.
\end{prf}

\smnoind
{\bf Proof of Theorem~\ref{thm:spec-oper-smiley}.}  The proof is divided into two cases.

\smnoind
{\bf Case $s=1$}.  von Neumann's double commutant theorem tells us that every $A\in B(\KH)$ is $(1,1)$-type Smiley operator, that is, $C_1(C_1(A))\subseteq VN(A)$.

\smnoind
{\bf Case $s\geqslant 2m-1$}. We prove this case in two steps.

First, suppose that $A\in\mathscr{S}_0(\KH)$ which means that $A$ is a normal operator, equivalently, $N=0$ in the canonical decomposition in Definition~\ref{defn:spec-oper-fin}.

For any $X\in C_1(A)$, that is, $\ad_A(X)=0$, the classical Fuglede's theorem \cite[Chapter IV, Theorem 4.10]{Sol2018} shows that $\ad_{A^*}(X)=0$. 
Now we write $A=\Re A+i\Im A$, where the real part $\Re A=\frac{A+A^{\ast}}{2}$ and the imaginary part $\Im A=\frac{A-A^{\ast}}{2i}$ are both self-adjoint operators.
So $\ad_{\Re A}(X)=\ad_{\Im A}(X)=0$. Similarly, Setting $X_1=\Re X$ and $X_2=\Im X$, by Lemma~\ref{lem:im-real-commut}, we have $\ad_{\Re A}(X_i)=\ad_{\Im A}(X_i)=0\ (i=1,2)$. Hence, we have that, $\ad_{A}(X_i)=0 \ (i=1,2)$, and then $\ad^s_{A}(X_i)$=0 for all $X_i$. Now, for any $B\in C_s(C_s(A))$, it is apparent that $\ad^s_{X_i}(B)=0\ (i=1,2)$, and hence $\ad_{X_i}(B)=0\ (i=1,2)$ by Lemma~\ref{lem:s-1}. Consequently, it follows that $\ad_X(B)=0$. Hence $B\in C_1(C_1(A))$. 
The argument above implies that $C_s(C_s(A))\subseteq C_1(C_1(A))$ for any normal operator $A\in B(\KH)$ and any positive integer $s\geqslant 2$. This, combined with {\bf Case $s=1$},  implies that \emph{every normal operator $A\in B(\KH)$ is $(s,s)$-type Smiley operator for any integer $s\geqslant 1$}.

Next,  we suppose that $A$ is an operator in $\mathscr{S}_m(\KH)$ ($m\geqslant 1$) with the canonical decomposition $A=S+N$. For any positive integer $s\geqslant 2m+1$ and $B\in C_s(C_s(A))$, we have that, if $X$ lies in $C_s(S)$, i.e., $\ad_{S}^{s}(X)=0$, then by Lemma~\ref{lem:s-1} $\ad_{S}(X)=0$, i.e., $X\in C_1(S)$. 
Since $N$ is a nilpotent operator of order $m+1$, i.e., $N^{m+1}=0$, and so we obtain
$$\ad_{N}^{2m+1}(X)=\sum_{j=0}^{2m+1}(-1)^{2m+1-j} {2m+1 \choose j}N^{j}XN^{2m+1-j}=0,$$
 where all ${2m+1 \choose j}$  are the combination numbers (cf. \cite{Rot1936} or \cite[p. 776]{Rob1961}).
Since $\ad_{S}(X)=0$ in the argument above, we have $(\CR_A-\CL_A)(X)=(\CL_N-\CR_N)(X)$, and so
$$\ad_A^{2m+1}(X)=(\CR_A-\CL_A)^{2m+1}(X)=(\CL_N-\CR_N)^{2m+1}(X)=\ad_N^{2m+1}(X)=0.$$
Consequently, due to $s\geqslant 2m+1$, we have $\ad_{A}^{s}(X)=0$.
Hence $\ad_{X}^{s}(B)=0$ holds since $B$ is an element in $C_s(C_s(A))$.
So we have $B\in C_s(C_s(S))$ because of $X\in C_s(S)$. 
Since $S$ is a normal operator, the first step above shows that $B\in C_1(C_1(S))$. Hence $C_s(C_s(A))\subseteq C_1(C_1(S))$ for any positive integer $s\geqslant 2m+1$. 

Now, if a bounded operator $X$ belongs to $C_1(A)$, then $X$ commutes with every spectral resolution for $A$ (cf. \cite[Lemma 3, Section 3, Chapter XVI]{Dun-Sch1988} or \cite[Page 226]{Dun1958}).
Then $X\in C_1(S)$ since $A$ and $S$ have the same spectral resolution. Applying $B\in C_1(C_1(S))$, we have $B\in C_1(X)$. Hence it follows that $B\in C_1(C_1(A))$, which means that
$C_s(C_s(A))\subseteq C_1(C_1(S))\subseteq C_1(C_1(A))$ for any positive integer $s\geqslant 2m+1$. 
In other words, {\bf Case $s\geqslant 2m+1$} ($m\geqslant 1$) reduces to {\bf Case $s=1$}. So Theorem~\ref{thm:spec-oper-smiley} is obtained.

\begin{rem}\label{rem:smiley-cpt-sa}
Let $A$ be a bounded normal operator on $\KH$.
By Theorem~\ref{thm:spec-oper-smiley}, we know that every operator $B\in C_s(C_s(A))$ for some positive integer $s$ lies in $VN(A)$. In other words, the operator $B$ is determined by $A$. 
More precisely, by Borel functional calculus (cf. \cite[pp.72]{Mur1990}), we can further see that there exist $f\in B_{\infty}(\sigma(A))$ such that $B=f(A)$, where $B_{\infty}(\sigma(A))$ be the $C^{\ast}$-algebra of all bounded Borel measurable complex valued functions on $\sigma(A)$. 

In particular, if $A$ be a compact self-adjoint operator on $\KH$, then $A$ has the canonical spectral decomposition (also known as diagonalization) $A=\sum_{i=1}^{\infty}\lambda_i P_i$ (see \cite[Theorem 5.1, Chapter II]{Con1990}). Applying the classical Smiley's theorem to the finite dimensional range space of every $P_i$,  we can show that $B=\sum_{i=1}^{\infty} g_{i}(A)\chi_{\{\lambda_{i}\}}(A)$,  where $g_{i}$ be some polynomial and $\chi_{\{\lambda_{i}\}}$ is the characteristics function of the singleton set $\{\lambda_i\}$.
\end{rem}

\medskip
\section{An example of non-spectral operators}\label{sec:ex}
\smallskip

In this section, as an example of non-spectral Smiley operators, we prove that every infinite dimensional unilateral operator is a proper $(k,l)$-type Smiley operator, where $k\in\{1,2\}$ and $k\leqslant l$.

\begin{defn}\cite[Proposition 2.10, Chapter II]{Con1990}\label{defn:unilateral-oper}
Let $\ell^2(\BN^+)$ be the separable Hilbert space with an orthonormal basis $\{e_{n}\}_{n=1}^{\infty}$, which consists of all complex number sequences $(x_{1},x_{2},...)$ satisfying $\sum^{\infty}_{i=1}\abs{x_i}^2<+\infty$. The unilateral shift operator $A$ on $\ell^2(\BN^+)$ is the operator in $\ell^2(\BN^+)$ defined by $A(e_{n})=e_{n+1}$ for all $n\in\BN^+$, equivalently, $A(x_{1},x_{2},...)=(0,x_{1},x_{2},...)$ for any sequence $(x_{1},x_{2},...)\in\ell^2(\BN^+)$. Note that $\sigma(A)=\overline{\mathbb{D}}$, where $\overline{\mathbb{D}}$ is the closed unit circle in the complex plane (cf. \cite[Proposition 6.5, Chapter VII]{Con1990}). So the unilateral shift operator on $\ell^2(\BN^+)$ is obviously not quasi-nilpotent operator.
\end{defn}

\begin{prop}\label{prop:unil-oper-nonspec}\cite[Corollary of Problem 147]{Hal1982}
The unilateral shift operator on $\ell^2(\BN^+)$ does not have any non-trivial reducing subspace, and so is not a spectral operator.
\end{prop}

\begin{prop}\label{prop:ex}
The unilateral shift operator $A$ on $\ell^2(\BN^+)$ is a proper $(k,l)$-type Smiley operator, where $k\in\{1,2\}$ and $k\leqslant l$.
\end{prop}

\begin{prf}
Since both $C_k(C_l(A))\subseteq C_l(C_l(A))$ and $C_k(C_l(A))\subseteq C_k(C_k(A))$ hold for every element $A\in B(\KH)$, and any two positive integers $k,l$ satisfying $k\leqslant l$, it suffices to prove that $C_2(C_2(A))\subseteq {\rm Pol}(A)$.

Taking an arbitrary sequence $\xi=(\xi_{1},\xi_{2},...)^{T}\in\ell^{2}(\mathbb{N^+})$ and an infinite dimensional complex matrix $X=(x_{ij})_{i,j\in\BN^+}$ , we have
\begin{equation}\label{401}
\aligned
\ad_{A}(X)\xi&=(AX-XA)\xi\\
&=\left(0,\sum_{j=1}^{\infty}x_{1j}\xi_{j},\sum_{j=1}^{\infty}x_{2j}\xi_{j},...\right)^{T}
-\left(\sum_{j=1}^{\infty}x_{1,j+1}\xi_{j},\sum_{j=1}^{\infty}x_{2,j+1}\xi_{j},...\right)^{T}\\
&=\left(-\sum_{j=1}^{\infty}x_{1,j+1}\xi_{j},\sum_{j=1}^{\infty}(x_{1j}-x_{2,j+1})\xi_{j},...,\sum_{j=1}^{\infty}(x_{kj}-x_{k+1,j+1})\xi_{j},...\right)^{T}
\endaligned
\end{equation}
and
\begin{equation}\label{402}
\aligned
\ad_{A}^{2}(X)\xi&=[A(AX-XA)-(AX-XA)A]\xi\\
&=\left(0,-\sum_{j=1}^{\infty}x_{1,j+1}\xi_{j},\sum_{j=1}^{\infty}(x_{1j}-x_{2,j+1})\xi_{j},...,\sum_{j=1}^{\infty}(x_{kj}-x_{k+1,j+1})\xi_{j},...\right)^{T}\\
&-\left(-\sum_{j=1}^{\infty}x_{1,j+2}\xi_{j},\sum_{j=1}^{\infty}(x_{1,j+1}-x_{2,j+2})\xi_{j},...,\sum_{j=1}^{\infty}(x_{k-1,j+1}-x_{k,j+2})\xi_{j},...\right)^{T}\\
&=\Bigg(\sum_{j=1}^{\infty}x_{1,j+2}\xi_{j},\sum_{j=1}^{\infty}(-2x_{1,j+1}+x_{2,j+2})\xi_{j},\sum_{j=1}^{\infty}(x_{1j}-2x_{2,j+1}+x_{3,j+2})\xi_j,\\
&\sum_{j=1}^{\infty}(x_{2j}-2x_{3,j+1}+x_{4,j+2})\xi_{j},...,\sum_{j=1}^{\infty}(x_{k-2,j}-2x_{k-1,j+1}+x_{k,j+2})\xi_{j},...\Bigg)^{T}
\endaligned
\end{equation}
Letting  $\ad_{A}^{2}(X)\xi=0$, by the arbitrariness of $\xi$, we have that, every matrix $X$  in $C_{2}(A)$ has the form of

\noindent 
$\begin{pmatrix}
x_{11} & x_{12} & 0 & 0 & 0 &0&\cdots&0&\cdots\\
x_{21} & x_{22} & 2x_{12} & 0 & 0 &0 &\cdots&0&\cdots\\
x_{31} & x_{32} & 2x_{22}-x_{11} & 3x_{12} & 0 &0&\cdots&0&\cdots\\
x_{41} & x_{42} & 2x_{32}-x_{21} & 3x_{22}-2x_{11} & 4x_{12} &0&\cdots&0&\cdots\\
x_{51} & x_{52} & 2x_{42}-x_{31} & 3x_{32}-2x_{21} & 4x_{22}-3x_{11} & 5x_{12} &\cdots&0&\cdots\\
\vdots & \vdots & \vdots & \vdots & \vdots & \vdots &\vdots&\vdots&\vdots\\
x_{k,1} & x_{k,2} & 2x_{k-1,2}-x_{k-2,1} & 3x_{k-2,2}-2x_{k-3,1}&4x_{k-3,2}-3x_{k-4,1} &5x_{k-4,2}-4x_{k-5,1}   &\cdots&kx_{12}&\cdots\\
\vdots &  \vdots  &  \vdots &  \vdots &  \vdots &  \vdots  &  \vdots &\vdots&\vdots\\
\end{pmatrix}.$

In the following, we compute $C_{2}(C_{2}(A))$.
Since $$(X\xi)(1)=x_{11}\xi_{1}+x_{12}\xi_{2}$$ 
and
$$(X\xi)(k)=x_{k1}\xi_{1}+x_{k2}\xi_{2}+\sum_{j=1}^{k-2}[(j+1)x_{k-j,2}-jx_{k-j-1,1}]\xi_{j+2}+kx_{12}\xi_{k+1} \ (k\geqslant 2),$$
setting $(B\psi)(k)=\sum_{j=1}^{\infty}\beta_{kj}\psi_{j}$ for all $\psi=\{\psi_{j}\}_{j=1}^{\infty}$ in $\ell^{2}(\mathbb{N})$, we have
\begin{equation}\label{403}
(BX\xi)(k)=\sum_{l=1}^{\infty}\beta_{kl}\left\{x_{l1}\xi_{1}+x_{l2}\xi_{2}+\sum_{j=1}^{l-2}[(j+1)x_{l-j,2}-jx_{l-j-1,1}]\xi_{j+2}+lx_{12}\xi_{l+1}\right\},
\end{equation}
and
\begin{equation}\label{404}
\aligned
\ (XB\xi)(k)&=x_{k1}(B\xi)(1)+x_{k2}(B\xi)(2)+\sum_{l=1}^{k-2}[(l+1)x_{k-l,2}-lx_{k-l-1,1}](B\xi)(l+2)\\
&+kx_{12}(B\xi)(k+1) \\
&=x_{k1}\left(\sum_{j=1}^{\infty}\beta_{1j}\xi_{j}\right)+x_{k2}\left(\sum_{j=1}^{\infty}\beta_{2j}\xi_{j}\right)\\
&+\sum_{l=1}^{k-2}[(l+1)x_{k-l,2}-lx_{k-l-1,1}]\left(\sum_{j=1}^{\infty}\beta_{l+2,j}\xi_{j}\right)+kx_{12}\left(\sum_{j=1}^{\infty}\beta_{k+1,j}\xi_{j}\right) \\
&=\sum_{j=1}^{\infty}\Bigg\{x_{k1}\beta_{1j}+x_{k2}\beta_{2j}+\sum_{l=3}^{k}[(l-1)x_{k-l+2,2}-(l-2)x_{k-l+1,1}]\beta_{l,j}+kx_{12}\beta_{k+1,j}\Bigg\}\xi_{j}.
\endaligned
\end{equation}

Consequently, it follows that
\begin{equation}\label{407}
\aligned
(XB-BX)\xi(k) &=\Bigg\{x_{k1}\beta_{11}+x_{k2}\beta_{21}+\sum_{l=3}^{k}[(l-1)x_{k-l+2,2}-(l-2)x_{k-l+1,1}]\beta_{l,1}\\
&+kx_{12}\beta_{k+1,1}-\sum_{l=1}^{\infty}\beta_{kl}x_{l1}\Bigg\}\xi_{1}+\Bigg\{x_{k1}\beta_{12}+x_{k2}\beta_{22}\\
&+\sum_{l=3}^{k}[(l-1)x_{k-l+2,2}-(l-2)x_{k-l+1,1}]\beta_{l,2}+kx_{12}\beta_{k+1,2}-\sum_{l=1}^{\infty}\beta_{kl}x_{l2}\Bigg\}\xi_{2}\\
&+\sum_{j=3}^{\infty}\Bigg\{x_{k1}\beta_{1j}+x_{k2}\beta_{2j}+\sum_{l=3}^{k}[(l-1)x_{k-l+2,2}-(l-2)x_{k-l+1,1}]\beta_{l,j}\\
&+kx_{12}\beta_{k+1,j}-(j-1)\beta_{k,j-1}x_{l2}+\sum_{l=j}^{\infty}\beta_{kl}[(j-1)x_{l-j+2,2}-(j-2)x_{l-j+1,1}]\Bigg\}\xi_{j}
\endaligned
\end{equation}
Considering a special case when $x_{11}=1,x_{ij}=0(i\neq1,j\neq1)$, we have
\begin{equation}\label{408}
X\psi=(\psi_{1},0,-\psi_{3},-2\psi_{4},...,-(k-2)\psi_{k},...)^{T},
\end{equation}
and
\begin{equation}\label{409}
\aligned
\ &(XB-BX)\xi=\Bigg(\beta_{12}\xi_{2}+\sum_{j=3}^{\infty}(j-1)\beta_{1,j}\xi_{j},-\beta_{21}\xi_{1}\\
&+\sum_{j=3}^{\infty}(j-2)\beta_{2,j}\xi_{j},
-2\beta_{31}\xi_{1}-\beta_{32}\xi_{2},...,(1-k)\beta_{k1}\xi_{1}-(k-2)\beta_{k2}\xi_{2},...\Bigg)^{T}.
\endaligned
\end{equation}
Hence we have
\begin{equation}\label{410}
\aligned
\ &(X(XB-BX))\xi=\Bigg(\beta_{12}\xi_{2}+\sum_{j=3}^{\infty}(j-1)\beta_{1,j}\xi_{j},\\
&0,2\beta_{31}\xi_{1}+\beta_{32}\xi_{2},...,(1-k)(2-k)\beta_{k1}\xi_{1}+(k-2)^{2}\beta_{k2}\xi_{2},...\Bigg)^{T},
\endaligned
\end{equation}
and
\begin{equation}\label{411}
\aligned
\ &((XB-BX)X)\xi=\Bigg(\sum_{j=3}^{\infty}(j-1)(2-j)\beta_{1,j}\xi_{j},\\
&-\beta_{21}\xi_{1}-\sum_{j=3}^{\infty}(j-2)^{2}\beta_{2,j}\xi_{j},-2\beta_{31}\xi_{1},...,(1-k)\beta_{k1}\xi_{1},...\Bigg)^{T}.
\endaligned
\end{equation}

If $X(XB-BX)=(XB-BX)X$ holds, then, by \eqref{410} and \eqref{411}, we have 
$$\beta_{1k}=0(k\neq1),\ \beta_{2k}=0(k\neq2)\ ,\beta_{k1}=0(k\neq1),\ \beta_{k2}=0(k\neq2).$$

Next we consider another case when $x_{12}=1,x_{ij}=0(i\neq1,j\neq2)$, from the similar procedures as above, we can obtain 
$$\beta_{kj}=0(k\neq j),\ \beta_{k+1,k+1}-\beta_{k,k}=\beta_{k,k}-\beta_{k-1,k-1}(k\geqslant 2).$$ 

That is $B$ is an infinite dimensional  bounded diagonal operator with the matrix form ${\rm diag}(\beta_{11},\beta_{22},...)$ satisfying $\beta_{k+1,k+1}-\beta_{k,k}=\beta_{k,k}-\beta_{k-1,k-1}\ (k\geqslant 2)$. Since all diagonal entries in a bounded diagonal operators are uniformly bounded (see \cite[Problem 61]{Hal1982} or \cite[Exercise 8, Section II.1]{Con1990}), it is immediately known that 
$$\beta_{k+1,k+1}-\beta_{k,k}=\beta_{k,k}-\beta_{k-1,k-1}=0 \ (k\geqslant 2).$$ 
So we can set $\beta_{k,k}=\lambda$ where $\lambda$ is an arbitrary complex number, and then $B=\lambda\cdot \id_{\ell^2(\BN^+)}$, which says that $$C_{2}(C_{2}(A))=\mathbb{C}I\subseteq{\rm Pol}(A).$$
The proof is done now.
\end{prf}



\bigskip

\begin{thebibliography}{99}

\bibitem{Con1990} J.\ B.\ Conway,  {\it A course in functional analysis}, 2nd ed., Graduate Texts in Mathematics {\bf 96}, Springer-Verlag, New York,1990.

\bibitem{Dun1952} N.\ Dunford,  The reduction problem in spectral theory. {\it Proceedings of the International Congress of Mathematicians, Cambridge, Mass.}, {\bf 2}(1950),115--122.  {\it Amer. Math. Soc., Providence, R. I.}, 1952.

\bibitem{Dun1958} N.\ Dunford, A survey of the theory of spectral operators, {\it Bull. Amer. Math. Soc.}, {\bf 64}(1958), 217--275.

\bibitem{Dun-Sch1988} N.\ Dunford, J.\ T.\ Schwartz, {\it Linear operators. Part III. Spectral operators}, Reprint of the 1971 original, Wiley Classics Library,  A Wiley-Interscience Publication, John Wiley and Sons, Inc., New York, 1988.

\bibitem{Hal1982} P.\ R.\ Halmos, {\it A Hilbert space problem book}, 2nd ed., Graduate Texts in Mathematics {\bf 19}, Springer-Verlag, New York, 1982.

\bibitem{Mur1990} G.\ J.\ Murphy, {\it $C^{\ast}$-algebras and operator theory}, Academic Press, San Diego, 1990.

\bibitem{Pra1994} V.\ Prasolov, {\it Problems and Theorems in Linear Algebra}, Translations of Mathematical Monographs {\bf 134}, American Mathematical Society, Providence, RI, 1994.

\bibitem{Rob1961} D.\ W.\ Robinson, A Note on k-Commutative Matrices, {\it J. Math. Phys.}, {\bf 2}(1961), 776--777.

\bibitem{Rob1964} D.\ W.\ Robinson, On matrix commutators of higher order, \emph{Canad. J. math.}, {\bf 17}
(1965), 527--532.

\bibitem{Rob1976} D.\ W.\ Robinson, From Pebbles to Commutators, {\it BYU Studies Quarterly}, {\bf 16}(1)(1976), Article 11. Available at: https://scholarsarchive.byu.edu/byusq/vol16/iss1/11.

\bibitem{Rot1936} W.\ E.\ Roth, On k-commutative matrices. {\it Trans. Amer. Math. Soc.}, {\bf 39}(3)(1936), 483--495.

\bibitem{Smi1961} M.\ F.\ Smiley, Matrix commutators, {\it Canad. J. Math.}, {\bf 13}(1961), 353--355.

\bibitem{Sol2018} P.\ Soltan, {\it A primer on Hilbert space operators}, Compact Textbooks in Mathematics, Birkh\"auser/Springer, Cham,  2018.

\bibitem{Wer1954} J.\  Wermer, Commuting spectral measures on Hilbert space, \emph{Pacific J. math.}, {\bf 4}
(1954), 355--361.





\end{thebibliography}
\end{document}